# New Iterative Methods for Interpolation, Numerical Differentiation and Numerical Integration

**M. Ramesh Kumar**
Phone No: +91 9840913580
Email address: ramjan_80@yahoo.com
Home page url: http//ramjan07.page.tl/

**Abstract**

Through introducing a new iterative formula for divided difference using Neville's and Aitken's algorithms, we study new iterative methods for interpolation, numerical differentiation and numerical integration formulas with arbitrary order of accuracy for evenly or unevenly spaced data. Basic computer algorithms for new methods are given.

## 1. Introduction

Interpolation is used in a wide variety of ways. Originally, it was used to do interpolation in tables defining common mathematical functions; but that is a far less important use in the present day, due to availability of computers and calculators. Interpolation is still use in the related problem of extending functions that are known only at a discrete set of points and such problems occurs frequently when numerically solving differential and integral equations. Next, Interpolation is used to solve problems from the more general area of interpolation theory. Interpolation is an important tool in producing computable approximations to commonly used functions. More over, to numerically integrate or differentiate a function, we often replace the function with simpler approximations and it is then integrated or differentiated. These simpler expressions are almost always obtained by interpolation [1].

A number of different methods have been developed to construct useful interpolation formulas for evenly or unevenly spaced points. Newton's divided difference formula [1,2,3], Lagrange's formula [1,2,3,10], Neville's and Aitken's iterated interpolation formulas[11,12] are the most popular interpolation formulas for polynomial interpolation to any arbitrary degree with finite number of points. The Lagrange formula is well suited for many theoretical uses of interpolation, but it is less desirable when actually computing the value of an interpolating polynomial [3]. Computation by Lagrange's method is quite laborious. For any computation the whole data is taken into calculation. If a new node is added, the computation has to be done afresh. These make Lagrange's method less suitable from the practical point of view. In this case, Neville's and Aitken's algorithms are very useful to iterate interpolation formula when a new node is added. Numerical computations by this method are simpler and less laborious than Lagrange's method. Also, it have an advantage over the Newton's interpolation formula is being very easily programmed for computer.

Numerical approximations to derivatives are used mainly in two ways. First, we are interested in calculating derivatives of given data that are often obtained empirically. Second, numerical differentiation formulae are used in deriving numerical methods for solving ordinary and partial differential equations [1]. A number of different methods have been developed to construct useful formulas for numerical derivatives. Most popular of the techniques are finite difference type [10], polynomial interpolation type [1,2,3,4], method of undetermined coefficients [1,2,3,8], and Richardson extrapolation [4,10]. More over, calculations of weights in finite difference formulas using recursive relations [7], explicit finite difference formulas [9] and few central difference formulas for finite and infinite data [5,6] are developed to construct useful numerical differentiation formulas. However, when a new node is added, the computation has to be done afresh. Thus, Iterative formula for numerical differentiation is still to be developed.

In Section 2, a new iterative formula for divided difference is for evenly or unevenly spaced data. In Section 3, the new iterative interpolation formulas are presented for both evenly or unequally spaced data with new divided difference table. In Section 4, new iterative methods for higher order numerical differentiation formulas are presented in recursive approach and also in direct form to any arbitrary order of accuracy for equally or unequally spaced data. In Section 5, the new iterated numerical integration formulas are derived from differentiation formulas, presented in Section 4 and the Taylor formula. In Section 6, is devoted to a brief conclusion.

## 2. Iterative Formula for divided difference

**Definition 2.1**

The $r^{\text{th}}$ divided difference of polynomial function $P(x)$ at the points $x, x_0, x_1, \ldots x_{r-1}$ is a polynomial in $x$, so we call it as divided difference polynomial of order r of $P(x)$. It is denoted by $P[x, x_0, x_1, \ldots x_{r-1}]$.

Now, Let

$$D_{j, j+1, \ldots, j+i}[x] = P[x, x_0, x_1, \ldots x_{r-1} / x_j, \ x_{j+1}, \ldots, x_{i+j}] \tag{2.1}$$

For $j = r, r+1, \ldots, n-i$ and $i = 0,1,2, \ldots, n-r$

Equation (2.1) is a divided difference polynomial iterated by the points $x_j, x_{j+1}, x_{j+2}, \ldots, x_{i+j}$

$$d_{r, r+1, \ldots i-1, i, j}[x] = P[x, x_0, x_1, \ldots x_{r-1} / x_r, \ x_{r+1}, \ldots, x_{i-1}, x_i, x_j] \tag{2.2}$$

For $j = i+1, i+2, \ldots, n$ and $i = r, r+1, \ldots, n-1$, Where $D_r[x] = d_r[x] = P[x, x_0, x_1, \ldots x_{r-1} / x_r]$

Equation (2.2) is a divided difference polynomial iterated by the points $x_r, x_{r+1}, x_{r+2}, \ldots, x_{i-1}, x_i, x_j$

**Note 2.1.** As a notation, $\{a, b, c, \ldots\}$ denotes the smallest interval containing all of real numbers $a, b, c, \ldots$ [1].

**Theorem 2.1.** *Let* $x_0, x_1, \ldots x_{r-1}$ *are 'r' numbers and* $x_r, x_{r+1}, \ldots, x_n$ *are* $(n-r+1)$ *distinct numbers in the interval* $[p, q]$, $r \le n$, $x \in [p, q]$ *and* $f \in C^{n+1}[p, q]$, *then* $\xi \in \{x, x_0, x_1, \ldots x_n\}$

$$f[x, x_0, x_1, \ldots, x_{r-1}] = D_{r, r+1, r+2, \ldots, n}[x] + \frac{f^{(n+1)}(\xi)}{(n+1)!} \prod_{i=r}^{n} (x - x_i) \tag{2.3}$$

$$f[x, x_0, x_1, \ldots, x_{r-1}] = d_{r, r+1, r+2, \ldots, n}[x] + \frac{f^{(n+1)}(\xi)}{(n+1)!} \prod_{i=r}^{n} (x - x_i) \tag{2.4}$$

**Proof.**

Let $P_n(x)$ is a polynomial of degree $\le n$ in '$x$' that approximates the function $f$ and takes the functional values $f(x_0), f(x_1), \ldots f(x_n)$ for the arguments $x_0, x_1, \ldots x_n$ respectively. It can be written as

$$P_n(x) = a_0 + a_1 x + a_2 x^2 + \ldots + a_n x^n, \quad \text{All } a's \in R \tag{2.5}$$

Then, we can write $r^{th}$ order divided difference of $P_n(x)$ at the points $x_0, x_1, \ldots x_{r-1}$ in terms of '$x$' as in the following form

$$P_n[x, x_0, x_1, \ldots x_{r-1}] = \hat{a}_0 + \hat{a}_1 x + \hat{a}_2 x^2 + \ldots + \hat{a}_{n-r} x^{n-r} = \tilde{P}_{n-r}(x) \text{ (Say)}, \quad \text{All } \hat{a}'s \in R \tag{2.6}$$

Now $\tilde{P}_{n-r}(x)$ is a polynomial of degree $\le n-r$. To interpolate it by Neville's method of iterated interpolation, we can use remaining points $x_{r,} x_{r+1}, \dots, x_n$ .

$$D_{j,j+1,\dots,j+i}[x] = \frac{1}{x_{i+j} - x_j} \begin{vmatrix} D_{j,j+1,\dots,j+i-1}[x] & x_j - x \\ D_{j+1,j+2,\dots,j+i}[x] & x_{i+j} - x \end{vmatrix}, \ j = r, r+1, \dots, n-i \text{ and } i = 0,1,2\dots, n-r \tag{2.7}$$

After $n-r$ iterations, we get $(n-r)^{th}$ degree polynomial $D_{r,r+1,\dots,n}[x]$.

i.e, $$P_n[x_i, x_0, x_1, x_2 \dots x_{r-1}] = D_{r,r+1,\dots,n}[x] \tag{2.8}$$

Similarly, by Aitken's method of iterated interpolation,

$$d_{r,r+1,\dots,i-1,i,j}[x] = \frac{1}{x_j - x_i} \begin{vmatrix} d_{r,r+1,\dots,i-1,i}[x] & x_i - x \\ d_{r,r+1,\dots,i-1,j}[x] & x_j - x \end{vmatrix}, \ \ j = i+1, i+2, \dots, n \ \text{ and } \ i = r, r+1, \dots, n-1 \tag{2.9}$$

After $n-r$ iterations, we get another $(n-r)^{th}$ degree polynomial $d_{r,r+1,\dots,n}[x]$.

i.e, $$P_n[x_i, x_0, x_1, x_2 \dots x_{r-1}] = d_{r,r+1,\dots,n}[x] \tag{2.10}$$

The above equation is polynomial form of $r^{th}$ divided difference of a polynomial. Since $P_n(x)$ approximates the function $f(x)$ on [p, q], so we have the following two equations (2.11) and (2.12). For some $\xi \in \{x, x_0, x_1, \dots x_n\}$

$$f(x) = P_n(x) + \frac{f^{(n+1)}(\xi)}{(n+1)!} \prod_{i=0}^{n} (x - x_i) \tag{2.11}$$

$$f[x_i, x_0, x_1, x_2 \dots, x_{r-1}] = P_n[x_i, x_0, x_1, x_2 \dots x_{r-1}], \quad i = r, r+1, r+2, , \dots, n \tag{2.12}$$

Using (2.11), we can write

$$f[x, x_0, x_1, \dots, x_{r-1}] = P_n[x, x_0, x_1, \dots, x_{r-1}] + \frac{f^{(n+1)}(\xi)}{(n+1)!} \prod_{i=r}^{n} (x - x_i) \tag{2.13}$$

Using (2.8) in (2.8), we get (2.3) and (2.10) in (2.13), we get (2.4)

## 3. Interpolation formulas.

**Theorem 3.1.** *Let $x_0, x_1, \dots x_{r-1}$ are 'r' numbers and $x_{r,} x_{r+1}, \dots, x_n$ are $(n-r+1)$ distinct numbers in the interval $[p,q]$, $r \le n$, $x \in [p,q]$ and $f \in C^{n+1}[p,q]$, for some $\xi \in \{x, x_0, x_1, \dots x_n\}$ then*

(i). $$f(x) = \sum_{i=0}^{r-1} f[x_0, x_1, \dots, x_i] \prod_{j=0}^{i-1} (x - x_j) + D_{r,r+1,r+2,\dots,n}[x] \prod_{i=0}^{r-1} (x - x_i) + \frac{f^{(n+1)}(\xi)}{(n+1)!} \prod_{i=0}^{n} (x - x_i) \tag{3.1}$$

(ii). $$f(x) = \sum_{i=0}^{r-1} f[x_0, x_1, \dots, x_i] \prod_{j=0}^{i-1} (x - x_j) + d_{r,r+1,r+2,\dots,n}[x] \prod_{i=0}^{r-1} (x - x_i) + \frac{f^{(n+1)}(\xi)}{(n+1)!} \prod_{i=0}^{n} (x - x_i) \tag{3.2}$$

**Proof.**

We know that

$$f[x, x_0, x_1, \ldots, x_{r-1}] = \frac{f[x, x_0, x_1, \ldots, x_{r-2}] - f[x_0, x_1, x_2 \ldots, x_{r-1}]}{x - x_{r-1}} \tag{3.3}$$

Rearranging this, we get

$$f[x, x_0, x_1, \ldots, x_{r-2}] = f[x_0, x_1, x_2 \ldots, x_{r-1}] + (x - x_{r-1}) f[x, x_0, x_1, \ldots, x_{r-1}] \tag{3.4}$$

Repeating this, we get

$$\begin{aligned} f(x) = f(x_0) &+ (x - x_0) f[x_0, x_1] + (x - x_0)(x - x_1) f[x_0, x_1, x_2] + \ldots \\ &+ (x - x_0)(x - x_1)(x - x_2) \ldots (x - x_{r-2}) f[x_0, x_1, x_2 \ldots, x_{r-1}] \\ &+ (x - x_0)(x - x_1) \ldots (x - x_{r-1}) f[x, x_0, x_1, \ldots, x_{r-1}] \end{aligned} \tag{3.5}$$

Using (2.3) in (3.5) and after simplification, we get (3.1), similarly, using (2.4) in (3.5) we get (3.2).

### 3.1. New Divided Difference Table with Neville's and Aitken's scheme

In Newton divided difference table, divided differences of new entries in each column are determined by divided difference of two neighboring entries in the previous column. But, the procedure of new divided difference table is different from the Newton divided difference table. For example, consider the argument values $x_0, x_1, x_2 \ldots, x_6$ for the corresponding functional values $f_0, f_1, f_2, \ldots, f_6$. As a matter of convenience, we write $f_k = f(x_k)$. The table 1 is divided into two parts. The first part of the table follows the procedure of New divided difference table and the second part of the table follows Neville scheme. The first order divided differences in the third column of the Table 1 are found by the sequence of evaluating $f[x_0, x_1]$, $f[x_0, x_2], \ldots,$ The second order divided differences in the fourth column of the Table 1 are found by the sequence of evaluating $f[x_0, x_1, x_2]$, $f[x_0, x_1, x_3], \ldots$. Similarly, the sequences $f[x_0, x_1, x_2, x_3]$, $f[x_0, x_1, x_2, x_4], \ldots$ are evaluated for fifth column. Then from the sixth column we use Neville's scheme.

Table 1. New divided difference table with Neville's Scheme

| **x** | **y** | $\bar{\delta}^1$ | $\bar{\delta}^2$ | $\bar{\delta}^3$ | Neville's Scheme to construct divided difference polynomial | | |
|---|---|---|---|---|---|---|---|
| $x_0$ | $f_0$ | | | | | | |
| | | $f[x_0, x_1]$ | | | | | |
| $x_1$ | $f_1$ | | $f[x_0, x_1, x_2]$ | | | | |
| | | $f[x_0, x_2]$ | | $f[x_0, x_1, x_2, x_3]$ | | | |
| $x_2$ | $f_2$ | | $f[x_0, x_1, x_3]$ | | $D_{3,4}[x]$ | | |
| | | $f[x_0, x_3]$ | | $f[x_0, x_1, x_2, x_4]$ | | $D_{3,4,5}[x]$ | |
| $x_3$ | $f_3$ | | $f[x_0, x_1, x_4]$ | | $D_{4,5}[x]$ | | $D_{3,4,5,6}[x]$ |
| | | $f[x_0, x_4]$ | | $f[x_0, x_1, x_2, x_5]$ | | $D_{4,5,6}[x]$ | |
| $x_4$ | $f_4$ | | $f[x_0, x_1, x_5]$ | | $D_{5,6}[x]$ | | |
| | | $f[x_0, x_5]$ | | $f[x_0, x_1, x_2, x_6]$ | | | |
| $x_5$ | $f_5$ | | $f[x_0, x_1, x_6]$ | | | | |
| | | $f[x_0, x_6]$ | | | | | |
| $x_6$ | $f_6$ | | | | | | |

Similar procedure for New divided difference with Aitken's Scheme

## 4. Formulas of Numerical Differentiation

**Definition 4.1.** Define

$$N_j^{(r)}[x] = \frac{1}{(x_j - x)^r}, \quad N_{j,j+1,\ldots,i+j}^{(r)}[x] = \frac{1}{x_{i+j} - x_j}\begin{vmatrix} N_{j,j+1,\ldots,i+j-1}^{(r)}[x] & x_j - x \\ N_{j+1,j+2,\ldots,i+j}^{(r)}[x] & x_{i+j} - x \end{vmatrix} \tag{4.1}$$

$$\tilde{N}_j^{(k)}[x] = \frac{f(x_j)}{(x_j - x)^k}, \quad \tilde{N}_{j,j+1,\ldots,i+j}^{(k)}[x] = \frac{1}{x_{i+j} - x_j}\begin{vmatrix} \tilde{N}_{j,j+1,\ldots,i+j-1}^{(k)}[x] & x_j - x \\ \tilde{N}_{j+1,j+2,\ldots,i+j}^{(k)}[x] & x_{i+j} - x \end{vmatrix} \tag{4.2}$$

$j = 0,1,\ldots,n-i$ and $i = 1,2,\ldots,n$

$N_{j,j+1,\ldots,i+j}^{(r)}[x]$ and $\tilde{N}_{j,j+1,\ldots,i+j}^{(k)}[x]$ are constructed by Neville's Algorithm for $r = 0,1,\ldots,k$

At the $n^{th}$ iteration, Let $\quad N_{0,1,\ldots,n}^{(r)}[x] = N_r(x)$, (4.3)

**Definition 4.2.** Define

$$A_{0,j}^{(r)}[x] = \frac{1}{x_j - x_0}\begin{vmatrix} \dfrac{1}{(x_0 - x)} & x_0 - x \\ \dfrac{1}{(x_j - x)^r} & x_j - x \end{vmatrix}, \quad A_{0,1,\ldots,i-1,i,j}^{(r)}[x] = \frac{1}{x_j - x_i}\begin{vmatrix} A_{0,1,\ldots i-1,i}^{(r)}[x] & x_i - x \\ A_{0,1,\ldots,i-1,j}^{(r)}[x] & x_j - x \end{vmatrix}, \tag{4.4}$$

$$\tilde{A}_{0,j}^{(k)}[x] = \frac{1}{x_j - x_0}\begin{vmatrix} \dfrac{f(x_0)}{(x_0 - x)^k} & x_0 - x \\ \dfrac{f(x_j)}{(x_j - x)^k} & x_j - x \end{vmatrix}, \quad \tilde{A}_{0,1,\ldots,i-1,i,j}^{(k)}[x] = \frac{1}{x_j - x_i}\begin{vmatrix} \tilde{A}_{0,1,\ldots,i-1,i}^{(k)}[x] & x_i - x \\ \tilde{A}_{0,1,\ldots,i-1,j}^{(k)}[x] & x_j - x \end{vmatrix} \tag{4.5}$$

$j = i+1, i+2,\ldots,n$ and $i = 0,1,r,\ldots,n-1$

$A_{0,1,\ldots,i-1,i,j}^{(r)}[x]$ and $\tilde{A}_{0,1,\ldots,i-1,i,j}^{(k)}[x]$ are constructed by Aitken's Algorithm for $r = 0,1,\ldots,k$

At the $n^{th}$ iteration ,Let $A_{0,1,\ldots,n}^{(r)}[x] = A_r(x)$ (4.6)

**Theorem 4.1.** *Let* $x, x_0, x_1, \ldots x_n$ *are* $(n+1)$ *distinct numbers in the interval* $[p,q]$ , $k \in W$ *and* $f \in C^{n+k+1}[p,q]$ *for some* $\xi \in \{x, x_0, x_1, \ldots x_n\}$ *then*

(i). $$\frac{f^{(k)}(x)}{k!} + \frac{f^{(k-1)}(x)}{k-1!}N_1(x) + \ldots + \frac{f(x)}{0!}N_k(x) = \tilde{N}_{0,1,2,\ldots,n}^{(k)}[x] + \frac{f^{(n+k+1)}(\xi)}{(n+k+1)!}\prod_{i=0}^{n}(x - x_i) \qquad 4.7)$$

(ii). $$\frac{f^{(k)}(x)}{k!} + \frac{f^{(k-1)}(x)}{k-1!}A_1(x) + \ldots + \frac{f(x)}{0!}A_k(x) = \tilde{A}_{0,1,2,\ldots,n}^{(k)}[x] + \frac{f^{(k+n+1)}(\xi)}{(k+n+1)!}\prod_{i=0}^{n}(x - x_i) \tag{4.8}$$

**Proof.**

Let $P_n(x)$ is a polynomial of degree $\le n$ in '$x$' that approximates the function $f$. Using Equation (2.3)

$$f[\underbrace{x,\ldots,x}_{k+1\ times}] = D_{0,1,2,\ldots,n}[x] + \frac{f^{(n+k+1)}(\xi)}{(n+k+1)!}\prod_{i=0}^{n}(x-x_i)\,, \tag{4.9}$$

Where $P_n[\underbrace{x,\ldots,x}_{k+1\ times}/\, x_0, x_1, x_2, \ldots, x_n] = D_{0,1,2,\ldots,n}[x]$

We can write,

$$P_n[\underbrace{x,\ldots,x}_{k\ times}, x_i] = \frac{P_n[x_i, \underbrace{x,\ldots,x}_{k-1\ times}] - P_n[\underbrace{x,\ldots,x}_{k\ times}]}{x_i - x}$$

$$= \frac{P_n[x_i, \underbrace{x,\ldots,x}_{k-1\ times}]}{x_i - x} - \frac{P_n{}^{(k-1)}(x)}{k-1!\,(x_i - x)}$$

$$= -\frac{P_n{}^{(k-1)}(x)}{k-1!\,(x_i - x)} - \frac{P_n{}^{(k-2)}(x)}{k-2!\,(x_i - x)^2} + \frac{P_n[x_i, \underbrace{x,\ldots,x}_{k-2\ times}]}{(x_i - x)^2}$$

Preceding this, we get

$$= -\frac{P_n{}^{(k-1)}(x)}{k-1!\,(x_i - x)} - \frac{P_n{}^{(k-2)}(x)}{k-2!\,(x_i - x)^2} - \ldots - \frac{P_n(x)}{0!\,(x_i - x)^k} + \frac{P_n(x_i)}{0!\,(x_i - x)^k} \tag{4.10}$$

Applying Equation (2.3) for two consecutive data $x_j, x_{j+1}$

$$P_n[\underbrace{x,\ldots,x}_{k+1\ times}\big/ x_j, x_{j+1}] = \frac{1}{x_{j+1} - x_j}\begin{vmatrix} P_n[\underbrace{x,\ldots,x}_{k\ times}\big/ x_j] & x_j \\ P_n[\underbrace{x,\ldots,x}_{k\ times}\big/ x_{j+1}] & x_{j+1} \end{vmatrix} \tag{4.11}$$

Using (4.10) in (4.11), using properties of determinants, we get

$$= -\frac{P_n{}^{(k-1)}(x)}{k-1!} N^{(1)}_{j,j+1}[x] - \frac{P_n{}^{(k-2)}(x)}{k-2!} N^{(2)}_{j,j+1}[x] - \ldots - \frac{P_n(x)}{0!} N^{(k)}_{j,j+1}[x] + \tilde{N}^{(k)}_{j,j+1}[x] \tag{4.12}$$

Similarly, for three consecutive data $x_j, x_{j+1}, x_{j+2}$

$$\begin{aligned} &P_n[\underbrace{x,\ldots,x}_{k+1\ times}\big/ x_j, x_{j+1}, x_{j+2}] \\ &\quad = -\frac{P_n{}^{(k-1)}(x)}{k-1!} N^{(1)}_{j,j+1,j+2}[x] - \frac{P_n{}^{(k-2)}(x)}{k-2!} N^{(2)}_{j,j+1,j+2}[x] - \ldots - P_n(x) N^{(k)}_{j,j+1,j+2}[x] + \tilde{N}^{(k)}_{j,j+1,j+2}[x] \end{aligned} \tag{4.13}$$

Proceeding this, for some $i$

$$P_n[\underbrace{x,\ldots,x}_{k+1\,times}/x_j, x_{j+1}, x_{j+2},,\ldots, x_{j+i}]$$

$$= -\frac{P_n{}^{(k-1)}(x)}{k-1!} N^{(1)}_{j,j+1,\ldots,i+j}[x] - \frac{P_n{}^{(k-2)}(x)}{k-2!} N^{(2)}_{j,j+1,\ldots,i+j}[x] - \ldots - P_n(x) N^{(k)}_{j,j+1,\ldots,i+j}[x] + \tilde{N}^{(k)}_{j,j+1,\ldots,i+j}[x] \quad (4.14)$$

Putting $i = n$ and $j = 0$, (i.e at $n^{th}$ iteration)

$$P_n[\underbrace{x,\ldots,x}_{k+1\,times}/x_0, x_1, x_2,,\ldots, x_n]$$

$$= -\frac{P_n{}^{(k-1)}(x)}{k-1!} N^{(1)}_{0,1,\ldots,n}[x] - \frac{P_n{}^{(k-2)}(x)}{k-2!} N^{(2)}_{0,1,\ldots,n}[x] - \ldots - P_n(x) N^{(k)}_{0,1,\ldots,n}[x] + \tilde{N}^{(k)}_{0,1,\ldots,n}[x] \quad (4.15)$$

Using Equation (4.3) and (4.9) in (4.15) and $P_n$ approximates the function $f$.

$$f[\underbrace{x,\ldots,x}_{k+1\,times}] - \frac{f^{(n+k+1)}(\xi)}{(n+k+1)!}\prod_{i=0}^{n}(x-x_i)$$

$$= -\frac{f^{(k-1)}(x)}{k-1!} N_1(x) - \frac{f^{(k-2)}(x)}{k-2!} N_2(x) - \ldots - \frac{f(x)}{0!} N_k(x) + \tilde{N}^{(k)}_{0,1,\ldots,n}[x] \quad (4.16)$$

After simplification, we get (4.3).

Now, Using (2.4)

$$f[\underbrace{x,\ldots,x}_{k+1\,times}] = d_{0,1,2,\ldots,n}[x] + \frac{f^{(n+k+1)}(\xi)}{(n+k+1)!}\prod_{i=0}^{n}(x-x_i), \quad (4.17)$$

Where $P_n[\underbrace{x,\ldots,x}_{k+1\,times}/ x_0, x_1, x_{2,\ldots,} x_n] = d_{0,1,2,\ldots,n}[x]$

for some $i$ and $j$

$$P_n[\underbrace{x,\ldots,x}_{k+1\,times}/x_0, x_1, x_2,,\ldots, x_i, x_j]$$

$$= -\frac{P_n{}^{(k-1)}(x)}{k-1!} A^{(1)}_{0,1,\ldots,i,j}[x] - \frac{P_n{}^{(k-2)}(x)}{k-2!} A^{(2)}_{0,1,\ldots,i,j}[x] - \ldots - P_n(x) A^{(k)}_{0,1,\ldots,i,j}[x] + \tilde{A}^{(k)}_{0,1,\ldots,i,j}[x] \quad (4.18)$$

at $n^{th}$ iteration

$$P_n[\underbrace{x,\ldots,x}_{k+1\,times}/x_0, x_1, x_2,,\ldots, x_n]$$

$$= -\frac{P_n{}^{(k-1)}(x)}{k-1!} A^{(1)}_{0,1,\ldots,n}[x] - \frac{P_n{}^{(k-2)}(x)}{k-2!} A^{(2)}_{0,1,\ldots,n}[x] - \ldots - P_n(x) A^{(k)}_{0,1,\ldots,n}[x] + \tilde{A}^{(k)}_{0,1,\ldots,n}[x] \quad (4.19)$$

Using Equation (4.6) and (4.17) in (4.19), $P_n$ approximates the function $f$ and after simplification

$$f[\underbrace{x,\ldots,x}_{k+1\,times}]-\frac{f^{(n+k+1)}(\xi)}{(n+k+1)!}\prod_{i=0}^{n}(x-x_i) = -\frac{f^{(k-1)}(x)}{k-1!}A_1(x)-\frac{f^{(k-2)}(x)}{k-2!}A_2(x)-\ldots-\frac{f(x)}{0!}A_k(x)+\tilde{A}^{(k)}_{0,1,\ldots,n}[x] \quad (4.20)$$

After simplification, we get (4.4),

**Theorem 4.2.** *Let* $x, x_0, x_1, \ldots x_n$ *are distinct numbers in the interval* $[p,q]$, $t \in W$ *and* $f \in C^{n+t+1}[p,q]$, *then*

(i).
$$\frac{f^{(t)}(x)}{t!} = a_0\tilde{N}^{(t)}_{0,1,2,\ldots,n}[x]+a_1\tilde{N}^{(t-1)}_{0,1,2,\ldots,n}[x]+a_2\tilde{N}^{(t-2)}_{0,1,2,\ldots,n}[x]\ldots+a_t\tilde{N}^{(0)}_{0,1,2,\ldots,n}[x] +\prod_{i=0}^{n}(x-x_i)\left(a_0\frac{f^{(t+n+1)}(\xi_0)}{(t+n+1)!}+a_1\frac{f^{(t+n)}(\xi_1)}{(t+n)!}+\ldots+a_t\frac{f^{(n+1)}(\xi_t)}{(n+1)!}\right) \quad (4.21)$$

$a_0 = 1$, $a_k = -(a_0N_k + a_1N_{k-1} + a_2N_{k-2} + \ldots + a_{k-1}N_1)$, $k = 1,2,\ldots,t$ and $\xi_i \in \{x, x_0, x_1, \ldots x_n\}$, $i = 0,1,2,\ldots,t$

(ii).
$$\frac{f^{(t)}(x)}{t!} = a_0\tilde{A}^{(t)}_{0,1,2,\ldots,n}[x]+a_1\tilde{A}^{(t-1)}_{0,1,2,\ldots,n}[x]+a_2\tilde{A}^{(t-2)}_{0,1,2,\ldots,n}[x]\ldots+a_t\tilde{A}^{(0)}_{0,1,2,\ldots,n}[x] +\prod_{i=0}^{n}(x-x_i)\left(a_0\frac{f^{(t+n+1)}(\xi_0)}{(t+n+1)!}+a_1\frac{f^{(t+n)}(\xi_1)}{(t+n)!}+\ldots+a_t\frac{f^{(n+1)}(\xi_t)}{(n+1)!}\right) \quad (4.22)$$

$a_0 = 1$, $a_k = -(a_0A_k + a_1A_{k-1} + a_2A_{k-2} + \ldots + a_{k-1}A_1)$, $k = 1,2,\ldots,t$ and $\xi_i \in \{x, x_0, x_1, \ldots x_n\}$, $i = 0,1,2,\ldots,t$

**Proof:**

We can write for some $t = 0,1,2,\ldots$ using equation (4.3)

$$\frac{f^{(t)}(x)}{t!} = a_0\tilde{N}^{(0)}_{0,1,2,\ldots,n}[x]+a_1\tilde{N}^{(1)}_{0,1,2,\ldots,n}[x]+a_2\tilde{N}^{(1)}_{0,1,2,\ldots,n}[x]+\ldots+a_t\tilde{N}^{(k)}_{0,1,2,\ldots,n}[x] +\prod_{i=0}^{n}(x-x_i)\left(a_0\frac{f^{(t+n+1)}(\xi_0)}{(t+n+1)!}+a_1\frac{f^{(t+n)}(\xi_1)}{(t+n)!}+\ldots+a_t\frac{f^{(n+1)}(\xi_t)}{(n+1)!}\right) \quad (4.23)$$

with all unknown $a$'s, To these $a$'s , Put $f(x) = 1$ in (4.23), If $t = 0$, then $a_0 = 1$ and for $t = 1,2,3,4,\ldots$ we have

$$a_0N_t + a_1N_{t-1} + a_2N_{t-2} + \ldots + a_tN_0 = 0 \quad (4.24)$$

Using (4.23) and (4.24) , we obtain (4.21).

Similarly, we can write for some $t = 0,1,2,\ldots$ using equation (4.4)

$$\frac{f^{(t)}(x)}{t!} = a_0\tilde{A}^{(t)}_{0,1,2,\ldots,n}[x]+a_1\tilde{A}^{(t-1)}_{0,1,2,\ldots,n}[x]+a_2\tilde{A}^{(t-2)}_{0,1,2,\ldots,n}[x]+\ldots+a_t\tilde{A}^{(0)}_{0,1,2,\ldots,n}[x] +\prod_{i=0}^{n}(x-x_i)\left(a_0\frac{f^{(t+n+1)}(\xi_0)}{(t+n+1)!}+a_1\frac{f^{(t+n)}(\xi_1)}{(t+n)!}+\ldots+a_t\frac{f^{(n+1)}(\xi_t)}{(n+1)!}\right) \quad (4.25)$$

with all unknown $a$'s. To these $a$'s, Put $f(x)=1$ in (4.24), If $t=0$, then $a_0=1$ and for $t=1,2,3,4,\dots$ we have

$$a_0 A_t + a_1 A_{t-1} + a_2 A_{t-2} + \dots + a_t A_0 = 0 \tag{4.26}$$

Using (4.25) and (4.26) , we obtain (4.22).

**Algorithm 4.1.** For the unevenly spaced points $x_0, x_1, x_2, \dots, x_n$ and known the functional values $f(x_i)$ at $x_i$, $i=0,1,2,\dots,n$ then the steps to use ($n$+1) point formula to estimate $t^{th}$ derivative of $f(x)$ at $x$ are

Step1: For $k=1$ *to* $n$ do calculate $N_k(x)$ using (4.1) and (4.3)
Step 2: For $k=1$ *to* $n$ do calculate $\tilde{N}^{(k)}_{0,1,2,\dots,n}[x]$ using (4.2)
Step3: $a_0=1$, For $k=1$ *to* $n$ do $a_k=-(a_0 N_k + a_1 N_{k-1} + a_2 N_{k-2} + \dots + a_{k-1} N_1)$
Step 4: Use (4.21) to find $t^{th}$ derivative at '$x$'.

**Algorithm 4.2.** For the unevenly spaced points $x_0, x_1, x_2, \dots, x_n$ and known the functional values $f(x_i)$ at $x_i$, $i=0,1,2,\dots,n$ then the steps to use ($n$+1) point formula to estimate $t^{th}$ derivative of $f(x)$ at $x$ are

Step1: For $k=1$ *to* $n$ do calculate $A_k(x)$ using (4.4) and (4.6)
Step 2: For $k=1$ *to* $n$ do calculate $\tilde{A}^{(k)}_{0,1,2,\dots,n}[x]$ using (4.5)
Step3: $a_0=1$, For $k=1$ *to* $n$ do $a_k=-(a_0 A_k + a_1 A_{k-1} + a_2 A_{k-2} + \dots + a_{k-1} A_1)$
Step 4: Use (4.22) to find $t^{th}$ derivative at '$x$'.

**Note 4.1.**

To find $t^{th}$ order differentiation, algorithms (5.1) and (5.2), costs $2.5(t^2+t+1)n^2$ operations.

## 5. Formulas for Integration

**Theorem 5.1.** *Let* $x_0, x_1, x_2, \dots, x_n$ *are the distinct numbers in the interval* $[p,q]$, $x, x+h \in [p,q]$, $h \neq 0$ *and if* $f \in C^{2n+1}[p,q]$, *with*

(i).

$$\int_x^{x+h} f(x)dx = \tilde{N}^{(0)}_{0,1,2,\dots,n}[x]\gamma_0 + \tilde{N}^{(1)}_{0,1,2,\dots,n}[x]\gamma_1 + \tilde{N}^{(2)}_{0,1,2,\dots,n}[x]\gamma_2 + \dots + \tilde{N}^{(n)}_{0,1,2,\dots,n}[x]\gamma_n + \prod_{i=0}^{n}(x-x_i)\left(\frac{f^{(n+1)}(\xi_n)}{(n+1)!}\gamma_0 + \frac{f^{(n+2)}(\xi_{n-1})}{n+2!}\gamma_1 + \dots + \frac{f^{(2n+1)}(\xi_0)}{(2n+1)!}\gamma_n\right) + O(h^{n+2}) \tag{5.1}$$

Where $a_0=1$, $a_k=-(a_0 N_k + a_1 N_{k-1} + a_2 N_{k-2} + \dots + a_{k-1} N_1)$, $k=1,2,\dots,n$

(ii).

$$\int_x^{x+h} f(x)dx = \tilde{A}^{(0)}_{0,1,2,\dots,n}[x]\gamma_0 + \tilde{A}^{(1)}_{0,1,2,\dots,n}[x]\gamma_1 + \tilde{A}^{(2)}_{0,1,2,\dots,n}[x]\gamma_2 + \dots + \tilde{A}^{(n)}_{0,1,2,\dots,n}[x]\gamma_n + \prod_{i=0}^{n}(x-x_i)\left(\frac{f^{(n+1)}(\xi_n)}{(n+1)!}\gamma_0 + \frac{f^{(n+2)}(\xi_{n-1})}{n+2!}\gamma_1 + \dots + \frac{f^{(2n+1)}(\xi_0)}{(2n+1)!}\gamma_n\right) + O(h^{n+2}) \tag{5.2}$$

Where $a_0=1$, $a_k=-(a_0 A_k + a_1 A_{k-1} + a_2 A_{k-2} + \dots + a_{k-1} A_1)$, $k=1,2,\dots,n$

$\gamma_k = \frac{h^{k+1}}{k+1} + a_1 \frac{h^{k+2}}{k+2} + \cdots + a_{n-k} \frac{h^{n+1}}{n+1}$, $k = 0,1,2,\ldots n$ and $\xi_i \in \{x, x_0, x_1, \ldots x_n\}$, $i = 0,1,2,\ldots,n$

**Proof.**

Using Taylor series on integration

$$\int_x^{x+h} f(x)dx = f(x)h + f'(x)\frac{h^2}{2!} + f''(x)\frac{h^3}{3!} + \ldots + f^{(n)}(x)\frac{h^{n+1}}{(n+1)!} + O(h^{n+2}) \tag{5.3}$$

Using (4.21) in equation (5.3) and simplifying we obtain,

$$\begin{aligned}\int_x^{x+h} f(x)dx = {} & \tilde{N}^{(0)}_{0,1,2,\ldots,n}[x]h + \left(\tilde{N}^{(1)}_{0,1,2,\ldots,n}[x] + \tilde{N}^{(0)}_{0,1,2,\ldots,n}[x]a_1\right)\frac{h^2}{2} \\ & + \left(\tilde{N}^{(2)}_{0,1,2,\ldots,n}[x] + \tilde{N}^{(1)}_{0,1,2,\ldots,n}[x]a_1 + \tilde{N}^{(0)}_{0,1,2,\ldots,n}[x]a_2\right)\frac{h^3}{3} + \ldots \\ & + \left(\tilde{N}^{(n)}_{0,1,2,\ldots,n}[x] + \tilde{N}^{(n-1)}_{0,1,2,\ldots,n}[x]a_1 + \ldots + \tilde{N}^{(1)}_{0,1,2,\ldots,n}[x]a_{n-1} + \tilde{N}^{(0)}_{0,1,2,\ldots,n}[x]a_n\right)\frac{h^{n+1}}{n+1} \\ & + \prod_{i=0}^{n}(x - x_i)\left(\frac{f^{(n+1)}(\xi_n)}{(n+1)!}\gamma_0 + \frac{f^{(n+2)}(\xi_{n-1})}{n+2!}\gamma_1 + \ldots + \frac{f^{(2n+1)}(\xi_0)}{(2n+1)!}\gamma_n\right) + O(h^{n+2})\end{aligned} \tag{5.4}$$

Rearranging the above equation,

$$\begin{aligned} = {} & \tilde{N}^{(0)}_{0,1,2,\ldots,n}[x]\left(h + a_1\frac{h^2}{2} + a_2\frac{h^3}{3} + \cdots + a_n\frac{h^{n+1}}{n+1}\right) + \tilde{N}^{(1)}_{0,1,2,\ldots,n}[x]\left(\frac{h^2}{2} + a_1\frac{h^3}{3} + \cdots + a_{n-1}\frac{h^{n+1}}{n+1}\right) \\ & + \tilde{N}^{(2)}_{0,1,2,\ldots,n}[x]\left(\frac{h^3}{3} + a_1\frac{h^4}{4} + \cdots + a_{n-2}\frac{h^{n+1}}{n+1}\right) + \ldots + \tilde{N}^{(n)}_{0,1,2,\ldots,n}[x]\frac{h^{n+1}}{n+1} \\ & + \prod_{i=0}^{n}(x - x_i)\left(\frac{f^{(n+1)}(\xi_n)}{(n+1)!}\gamma_0 + \frac{f^{(n+2)}(\xi_{n-1})}{n+2!}\gamma_1 + \ldots + \frac{f^{(2n+1)}(\xi_0)}{(2n+1)!}\gamma_n\right) + O(h^{n+2})\end{aligned} \tag{5.5}$$

$$\begin{aligned} = {} & \tilde{N}^{(0)}_{0,1,2,\ldots,n}[x]\gamma_0 + \tilde{N}^{(1)}_{0,1,2,\ldots,n}[x]\gamma_1 + \tilde{N}^{(2)}_{0,1,2,\ldots,n}[x]\gamma_2 + \ldots + \tilde{N}^{(n)}_{0,1,2,\ldots,n}[x]\gamma_n \\ & \prod_{i=0}^{n}(x - x_i)\left(\frac{f^{(n+1)}(\xi_n)}{(n+1)!}\gamma_0 + \frac{f^{(n+2)}(\xi_{n-1})}{n+2!}\gamma_1 + \ldots + \frac{f^{(2n+1)}(\xi_0)}{(2n+1)!}\gamma_n\right) + O(h^{n+2})\end{aligned} \tag{5.6}$$

Thus, we obtain equation (5.1). Similarly using (4.22) in (5.3), after simplification, we get (5.2).

**Algorithm 5.1.** For the unevenly spaced points $x_0, x_1, x_2, \ldots, x_n$ and known the functional values $f(x_i)$ at $x_i$, $i = 0,1,2,\ldots,n$ then the steps to estimate $\int_x^{x+h} f(x)dx$ are,

Step1: For $k = 1$ *to* $n$ do calculate $N_k(x)$ using (4.1) and (4.3)

Step 2: For $k = 1$ *to* $n$ do calculate $\tilde{N}^{(k)}_{0,1,2,\ldots,n}[x]$ using (4.2)

Step3: $a_0 = 1$, For $k = 1$ *to* $n$ do $a_k = -(a_0N_k + a_1N_{k-1} + a_2N_{k-2} + \ldots + a_{k-1}N_1)$

Step 4: For $k = 0$ *to* $n$ do $\gamma_k = \frac{h^{k+1}}{k+1} + a_1 \frac{h^{k+2}}{k+2} + \cdots + a_{n-k} \frac{h^{n+1}}{n+1}$

Step 5: Use (5.1), to find numerically $\int_x^{x+h} f(x)dx$

**Algorithm 5.2.** For the unevenly spaced points $x_0, x_1, x_2, \ldots, x_n$ and known the functional values $f(x_i)$ at $x_i$, $i = 0,1,2,\ldots,n$ then the steps to estimate $\int_x^{x+h} f(x)dx$ are,

Step1: For $k = 1$ *to* $n$ do calculate $A_k(x)$ using (4.4) and (4.6)

Step 2: For $k = 1$ *to* $n$ do calculate $\tilde{A}^{(k)}_{0,1,2,\ldots,n}[x]$ using (4.5)

Step3: $a_0 = 1$, For $k = 1$ *to* $n$ do $a_k = -\left(a_0 A_k + a_1 A_{k-1} + a_2 A_{k-2} + \ldots + a_{k-1} A_1\right)$

Step 4: For $k = 0$ *to* $n$ do $\gamma_k = \frac{h^{k+1}}{k+1} + a_1 \frac{h^{k+2}}{k+2} + \cdots + a_{n-k} \frac{h^{n+1}}{n+1}$

Step 5: Use (5.2), to find numerically $\int_x^{x+h} f(x)dx$

## 6. Conclusion

By introducing a new iterated method for divided difference and new divided difference table, we have studied iterated methods for interpolation, numerical differentiation and integration formulas with arbitrary order accuracy for evenly or unevenly spaced data using Neville's and Aitken's algorithms. First, we study iterated interpolation formula which generalizes Newton interpolation formula and Iterated interpolation formula. However when a new node is added, we have to add one more data to new divided difference table. But new iterated formulas for higher order derivatives and numerical integration to arbitrary order of accuracy are very handier even when we add new data for further iteration. Basic computer algorithms are given for new formulas. Through new iterative method for divided difference, we have studied three major problems of Numerical analysis.

## List of References